\documentclass[graybox]{svmult}

\usepackage{type1cm}        % activate if the above 3 fonts are
\usepackage{makeidx}         % allows index generation
\usepackage{graphicx}        % standard LaTeX graphics tool
\usepackage{multicol}        % used for the two-column index
\usepackage[bottom]{footmisc}% places footnotes at page bottom

\usepackage{amssymb}
\usepackage{amsmath, amsthm}
\usepackage{subfigure}

\usepackage{newtxtext}       % 
\usepackage[varvw]{newtxmath}       % selects Times Roman as basic font

\newcommand{\RR}{{\mathbb{R}}}

\makeindex             % used for the subject index
\begin{document}

\title*{Using Multiple Dirac Delta Points to Describe Inhomogeneous Flux Density over a Cell Boundary in a Single-Cell Diffusion Model}
\titlerunning{Multiple Dirac Delta Points in a Diffusion Model}
% Use \titlerunning{Short Title} for an abbreviated version of
% your contribution title if the original one is too long
\author{Qiyao Peng\orcidID{0000-0002-7077-0727} and\\ Sander C. Hille\orcidID{0000-0003-0437-6745}}
% Use \authorrunning{Short Title} for an abbreviated version of
% your contribution title if the original one is too long
\institute{Qiyao Peng and Sander C. Hille \at Mathematical Institute, Leiden University.  Gorlaeus Building, Einsteinweg 55, 2333 CC, Leiden, The Netherlands. \email{q.peng@math.leidenuniv.nl}}
% \and Name of Second Author \at Name, Address of Institute \email{name@email.address}}
%
% Use the package "url.sty" to avoid
% problems with special characters
% used in your e-mail or web address
%
\maketitle

\abstract{Biological cells can release compounds into their direct environment, generally inhomogeneously over their cell membrane, after which the compounds spread by diffusion. In mathematical modelling and simulation of a collective of such cells, it is theoretically and numerically advantageous to replace spatial extended cells with point sources, in particular when cell numbers are large, but still so small that a continuum density description cannot be justified, or when cells are moving. We show that inhomogeneous flux density over the cell boundary may be realized in a point source approach, thus maintaining computational efficiency, by utilizing multiple, clustered point sources (and sinks). In this report, we limit ourselves to a sinusoidal function as flux density in the spatial exclusion model, and we show how to determine the amplitudes of the Dirac delta points in the point source model, such that the deviation between the point source model and the spatial exclusion model is small.}

\section{Introduction}
\label{sec:intro}
Cell-to-cell communication plays an important role in multi-cellular organisms, for example during development \cite{alberts2002molecular}, in attracting a proper immune system response (by secretion of so-called cytokines, that attract immune cells \cite{Peng_2020}, or by ensuring that cellular events proceed as expected \cite{Perbal_2003}). There are various methods for communication. One of them is via signaling molecules (like the cytokines), which are secreted, diffuse in the environment and are detected by receptor proteins on another cell's membrane. Receptors and other proteins are \textit{not} uniformly distributed over the cell membrane \cite{cooper2022cell}. As a result, the secretion of signaling molecules into the extracellular environment is expected to be inhomogeneous over the membrane.

One may view the intracellular environment and the extracellular environment as separated, spatially extended subdomains, particularly when the aim is to depict the process on a small spatial scale. We call this a \textit{spatial exclusion model}, because molecules in the environment cannot pass simply by diffusion through the space occupied by the cell. Computationally, using a Finite Element Method (FEM), this view is tractable when there are few static cells in the domain. Movement requires remeshing in every time iteration in FEM, which is time-consuming. On the other hand, at a larger spatial scale, cells are often reduced to `point masses' and a Dirac delta point source replaces flux boundary conditions at the no-longer extant cell membrane. This is done to improve the computational efficiency and simplify the spatial exclusion model. In particular, it will allow effective computation of diffusive communication between moving cells. We call this approach the \textit{point source model}. Similar upscalling work has been done in \cite{Peng_Vermolen_2022} concerning the linear momentum balance equation. At an even larger scale one uses a density description for the cell population. This shall not be considered here.

The spatial exclusion and point source approach lead to a difference between the solutions to these two models of different scale. It seems a general assumption among researchers that the deviation is small. To the best of our knowledge, this difference was only investigated in detail, quantitatively, in \cite{HMEvers2015} in a theoretical manner and in our previous work \cite{Peng2023}. In particular, in \cite{Peng2023}, we focused on a homogeneous flux density over the cell membrane. We found a systemic time delay related to the cell size relative to the strength of diffusion. We investigated how to compensate for this time delay when the diffusion coefficient is small.

Here, we consider the question: how to replace a spatial exclusion model with {\it inhomogeneous flux density} over the boundary, by a point source model? By its very nature, a point source spreads (or takes) `mass' equally in all directions. We shall show that one may proceed by replacing the spatially extended cell by a small cluster of point sources and sinks, with well-tuned intensity. 
% \vskip 0.2cm

The manuscript is structured as below: Section \ref{sec:math formulation} provides the mathematical formulation of the spatial exclusion and point source model and a corresponding phrasing of the main question. Section \ref{sec:model} discusses how to locate the Dirac points and compute their intensity. Numerical results are displayed in Section \ref{sec:results}. Conclusions and discussions are shown in Section \ref{sec:conclusion}.

\section{Mathematical formulation of the question}
\label{sec:math formulation}
\begin{figure}
    \centering
    \includegraphics[width = 0.5\textwidth]{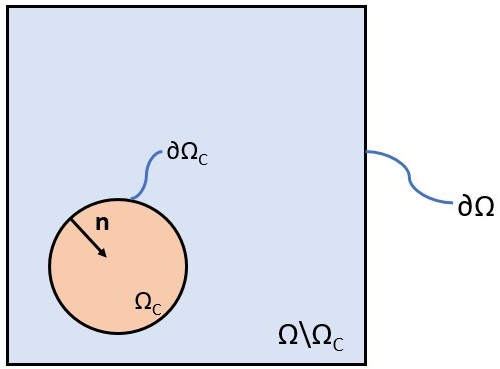}
    \caption{\it A schematic presentation of a single cell ($\Omega_C$ in orange) within domain ($\Omega$ as a square). $\Omega\setminus\Omega_C$ (blue domain) is the extracellular environment, into which the cell releases diffusive compounds that cannot escape through the boundary $\partial\Omega$. A typical direction of the normal vector $\boldsymbol{n}$ on $\partial\Omega_C$ is indicated.}
    \label{fig:domain_pic}
\end{figure}

The two types of model are mathematically formulated as follows, and Figure \ref{fig:domain_pic} depicts the domains of the models. The spatial exclusion model excludes the cell domains, $\Omega_{C_i}$, from the entire domain $\Omega$, which are both bounded and simply connected subsets of $\RR^2$ with piece-wise smooth boundary. Let $\Omega_C:=\bigcup_{i=1}^N \Omega_{C_i}$ be the total cell domain. At its boundary $\partial\Omega_C$, i.e. the total boundary of the cells, a flux boundary condition is prescribed. Explicitly:
\begin{equation}
	\label{Eq_BVP_hole}
	(BVP_S)\quad\left\{
	\begin{aligned}
	\frac{\partial u_S(\boldsymbol{x},t)}{\partial t} - D\Delta u_S(\boldsymbol{x},t) &= 0, &\mbox{in $\Omega\setminus\bar{\Omega}_C, t>0$,}\\
	D\nabla u_S(\boldsymbol{x},t)\cdot\boldsymbol{n} &= \phi(\boldsymbol{x},t), &\mbox{on $\partial\Omega_C, t>0$,}\\
	D\nabla u_S(\boldsymbol{x},t)\cdot\boldsymbol{n} &= 0, &\mbox{on $\partial\Omega, t>0$,}\\
	u_S(\boldsymbol{x}, 0) &= u_0(\boldsymbol{x}), &\mbox{in $\Omega\setminus\bar{\Omega}_C$,}
	\end{aligned}
	\right.
\end{equation}
where $D$ is diffusion coefficient of compounds, $u_S(\boldsymbol{x})$ is the concentration of the compounds (which is to be solved) defined by $(BVP_S)$, $u_0(\boldsymbol{x})$ is the initial concentration of the compounds over the computational domain, $\boldsymbol{n}$ is the outward pointing unit normal vector to the domain boundary of $\Omega\setminus\bar{\Omega}_C$. Note that the flux density $\phi(\boldsymbol{x},t)$ is non-negative at $\boldsymbol{x}\in\partial\Omega_C$ where there is a flux of compound into the environment of the cell, i.e. $\Omega\setminus\bar{\Omega}_C$. Furthermore, the initial-boundary value problem defined by point sources is then given by
\begin{equation}
	\label{Eq_BVP_dirac}
	(BVP_P)\quad\left\{
	\begin{aligned}
	\frac{\partial u_P(\boldsymbol{x},t)}{\partial t} - D\Delta u_P(\boldsymbol{x},t) &= \sum_{i = 1}^{N}\Phi_i(t)\delta(\boldsymbol{x}-\boldsymbol{x}_i), &\mbox{in $\Omega, t>0$,}\\
	D\nabla u_P\cdot\boldsymbol{n} &= 0, &\mbox{on $\partial\Omega, t>0$,}\\
	u_P(\boldsymbol{x}, 0) &= \bar{u}_0(\boldsymbol{x}), &\mbox{in $\Omega, t=0$.}
	\end{aligned}
	\right.
\end{equation}
Here, $u_P(\boldsymbol{x})$ is the concentration of compounds over the computational domain $\Omega$, $\bar{u}_0(\boldsymbol{x})$ is considered as a spatial extension of $u_0(\boldsymbol{x})$ in $(BVP_S)$ from $\Omega\setminus\Omega_C$ to $\Omega$, $\Phi_i(t)$ is a function that describes the flux of mass per unit time from the source at $\boldsymbol{x}_i$. Exact agreement of the solutions to these two models has been proven in \cite{Peng2023}, Proposition 1.1, when the cell boundary fluxes  in both agree:
\begin{equation}\label{PropCondition}
\phi(\boldsymbol{x},t) - D\nabla u_P(\boldsymbol{x},t)\cdot\boldsymbol{n}= 0, \qquad\mbox{ a.e. on $\partial\Omega_C\times[0,\infty)$.}
\end{equation}

In summary, to maintain good agreement between the two solutions in the extracellular environment, the deviation in flux over $\partial\Omega_C$ needs to be minimized. In the point source model the flux over the cell boundary is a result of diffusion from the source; it is `post-processed'. Hence, a natural question is, how to define the intensities of the Dirac delta points in the point source model such that the flux condition in Equation \eqref{PropCondition} can be best approximated.

\section{Point source localization and intensity}
\label{sec:model}
In this study, we mainly consider a prescribed inhomogeneous flux density distribution $\phi(\boldsymbol{x},t)$ over the circular boundary of a single cell centred at $\boldsymbol{x}_0$ of radius $R$ that is time-independent and in the form of 
\begin{equation}
   \label{Eq_phi}
   \phi(\boldsymbol{x}_\theta) = \phi(\theta) = \phi_0+A\sin(n\theta) = \phi_0(1+\rho\sin(n\theta)),
\end{equation}
in angular coordinate $\theta$. That is, $\boldsymbol{x}_\theta := \boldsymbol{x}_0 + R(\cos\theta,\sin\theta)\in\partial\Omega_C$. In this manuscript, we only consider a single Fourier mode. More general spatially inhomogeneous flux densities can be obtained in principle by superposition, using Fourier series expansion. Moreover, $\phi_0, A>0$, $\rho:=\frac{A}{\phi_0}\leqslant 1$ and $n$ a positive integer. We propose an approach to locate and compute the intensities of the Dirac delta points. Under certain conditions, the emergent boundary flux in the point source model converges to the predefined flux density of the spatial exclusion model.

The idea to achieve such an inhomogeneous flux is by putting one point source at the centre $\boldsymbol{x}_0$, essentially realising a homogenous flux on the boundary. The spatial variation is then realised by the non-centred Dirac delta points. The variation in $\phi(\boldsymbol{x})$ around $\phi_0$, i.e. $A\sin(n\theta)$, takes extreme values at $\theta_k:=(k-\frac{1}{2})\frac{\pi}{n}$, $k=1,2,\dots,2n$, with a maximum at odd $k$ and a minimum at even $k$.
For the sake of keeping the model simple, we want to have as small a number of point sources as possible, while retaining a good approximation of the spatial exclusion model with inhomogeneous flux $\phi$ over $\partial\Omega_C$.

Suppose that the additional $N$ non-centred Dirac delta point sources are located at $\boldsymbol{x}_i\in\Omega_C$, $i=1,\dots, N$. Let $\Phi_i(t)$ be the intensity of the point source (or sink, if $\Phi_i(t)<0$) at $\boldsymbol{x}_i$. If $\Omega=\RR^2$, or $\partial\Omega$ is remote from the cell and $t$ is small such that reflection of mass from the boundary $\partial\Omega$ can still be neglected, the diffusion system defined in \eqref{Eq_BVP_dirac}, yields a flux density over the cell boundary that is given by
\begin{equation}
    \label{Eq_bnd_flux_general_s}
    \begin{aligned}
    D\nabla u_P(\boldsymbol{x},t)\cdot\boldsymbol{n} &= \sum_{i = 0}^N \int_0^t \frac{\Phi_i(s)}{4\pi D(t-s)}\exp\left\{-\frac{\|\boldsymbol{x} - \boldsymbol{x}_i\|^2}{4D(t-s)}\right\}\frac{(\boldsymbol{x} - \boldsymbol{x}_i)\cdot(\boldsymbol{x}-\boldsymbol{x}_C)}{2(t-s)\|\boldsymbol{x}-\boldsymbol{x}_C\|} ds,\\
    \end{aligned}
\end{equation}
for any point $\boldsymbol{x}\in\partial\Omega_C$ and $\boldsymbol{x}_C$ presents the cell centre.

Generally, equality of $\phi(\boldsymbol{x})$ to the expression in \eqref{Eq_bnd_flux_general_s}, ensuring perfect agreement of solutions according to \eqref{PropCondition}, cannot hold
at all $\boldsymbol{x}\in\partial\Omega_C$ at all time. Not even if $A=0$, see \cite{Peng2023}. Thus, approximation of some sort is required. The convolution integrals in \eqref{Eq_bnd_flux_general_s} cause concern for the analysis. Therefore, for each $t$, we write
\begin{equation}\label{eq:Phi deviation mean}
    \Phi_i(s)\ =\ \widetilde{\Phi}_i(t) \ +\ \delta_t(s),\qquad s\in[0,t].
\end{equation}
That is, we view $\Phi_i(s)$ on $[0,t]$ as a deviation $\delta_t(s)$ from 'the mean' $\widetilde{\Phi}_i(t)$. Replacing $\Phi_i(s)$ by $\widetilde{\Phi}_i(t)$ and performing the integration in \eqref{Eq_bnd_flux_general_s} yields the expression
\begin{equation}\label{eq:expresion tilde Phi}
    \tilde{\phi}(\boldsymbol{x},t) :=  \sum_{i=0}^N \frac{\widetilde{\Phi}_i(t)}{2\pi R}  \frac{(\boldsymbol{x} - \boldsymbol{x}_C)\cdot(\boldsymbol{x}-\boldsymbol{x}_i)}{\|\boldsymbol{x} - \boldsymbol{x}_i\|^2}   \exp\left\{-\frac{\|\boldsymbol{x} - \boldsymbol{x}_i\|^2}{4Dt}\right\}
\end{equation}

We ignore the contribution of $\delta_t(s)$ to the point source flux $D\nabla u_P(\boldsymbol{x},t)\cdot\boldsymbol{n}$ in  \eqref{Eq_bnd_flux_general_s} and require the flux densities $\tilde{\phi}(\boldsymbol{x},t)$ to be equal to $\phi(\boldsymbol{x})$ for all time $t>0$ at the extreme points of $\phi$ only. That is,
\begin{equation}\label{eq:conditions for Phi}
    \tilde{\phi}(\boldsymbol{x}_{\theta_k},t) = \phi(\boldsymbol{x}_{\theta_k}) = \begin{cases} \phi_0 + A, & \mbox{if $k$ is odd},\\
    \phi_0 - A, & \mbox{if $k$ is even},
    \end{cases}
\end{equation}
for $k=1,2,\dots, 2n$. Thus one obtains $2n$ equations for $N+1$ unknown $\widetilde{\Phi}_i(t)$. So, with general locations $\boldsymbol{x}_i$ for the point sources in $\Omega_C$ one needs $N=2n-1$ points to be able to solve \eqref{eq:conditions for Phi}. However, by exploiting symmetry one may specifically localize point sources in such a way that one can do with less points.

\begin{figure}
    \centering
    \subfigure[$\phi(\theta) = \phi_0+A\sin(\theta)$]{
    \includegraphics[width = 0.32\textwidth]{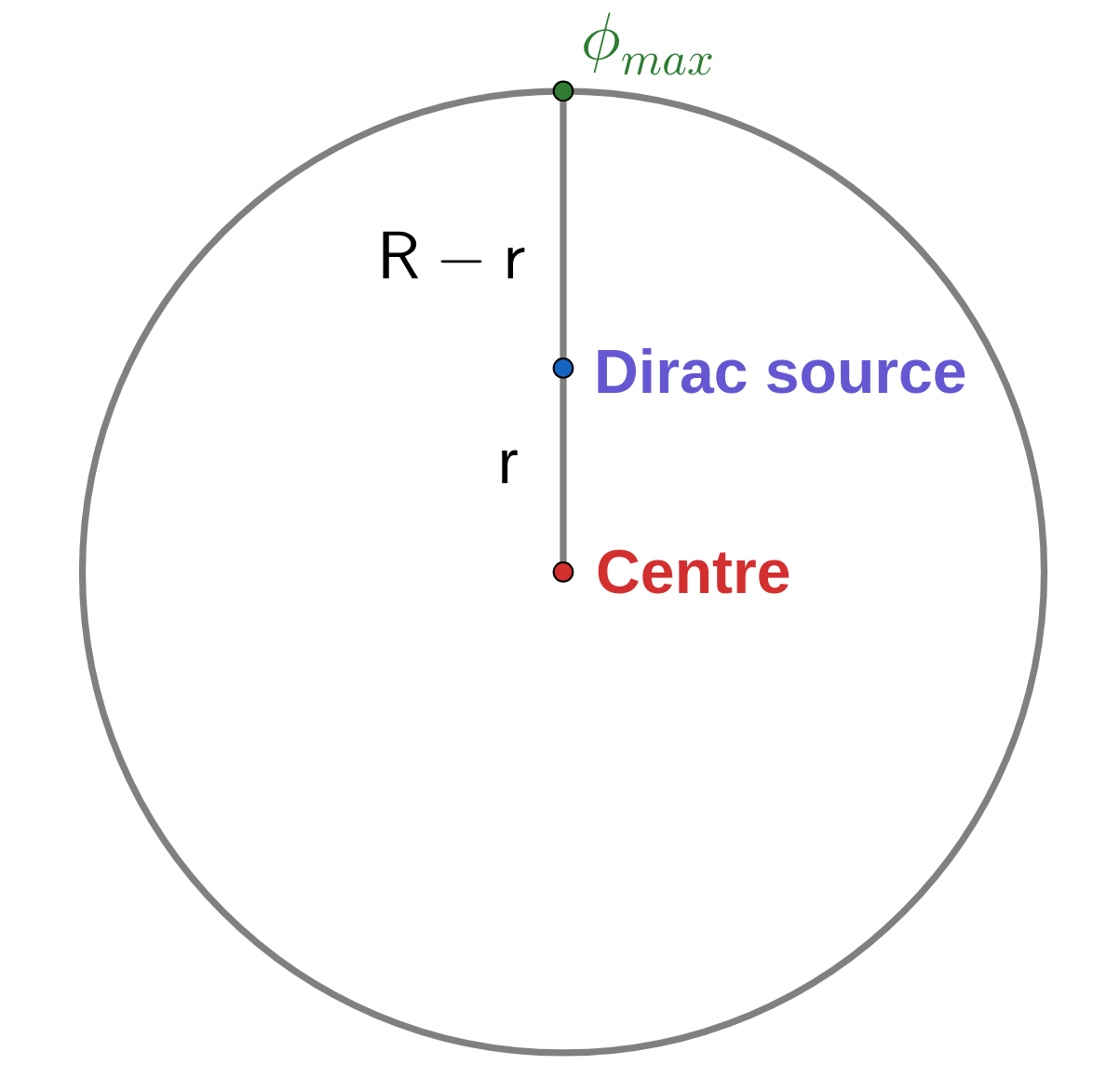}}
    \hfil
    \subfigure[$\phi(\theta) = \phi_0+A\sin(2\theta)$]{
    \includegraphics[width = 0.32\textwidth]{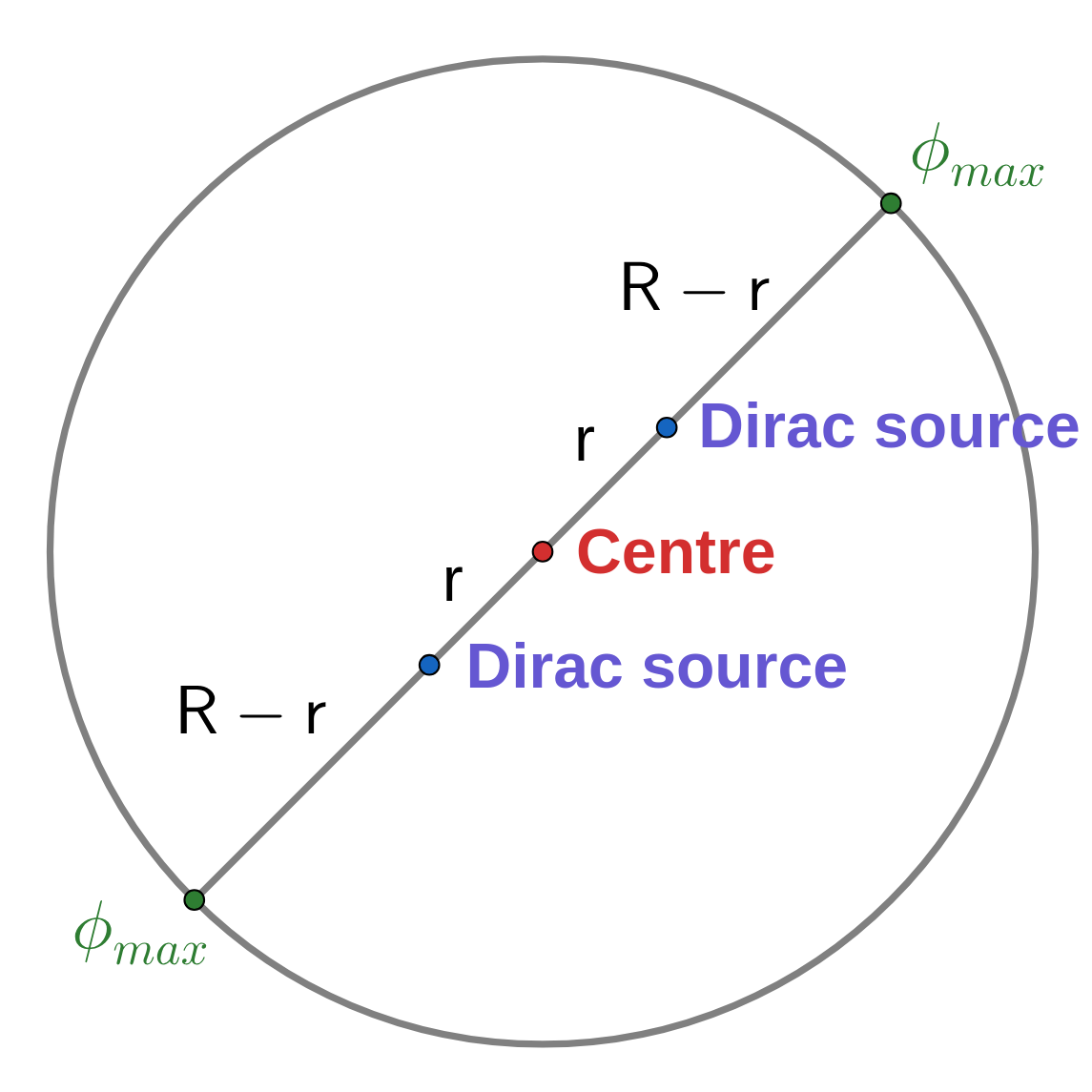}}
    \caption{{\it The locations of Dirac delta points for $n=1$ and $n=2$ are shown in Panel (a) and (b), respectively. The distance between the non-centred Dirac delta point and the center is denoted as $r$ and all the non-centred Dirac delta points are identical.}}
    \label{fig:1_2_Dirac}
\end{figure}

We locate the off-centre points each on one of the $n$ line segments connecting the centre $\boldsymbol{x}_0$ to the location $\boldsymbol{\theta_{k}}$ of the maximum of $\phi(\boldsymbol{x})$, for $k$ odd, at a distance $r$ from the centre.
See Figure \ref{fig:1_2_Dirac} for the examples of $n=1$ and $n=2$. Moreover, due to the symmetry of $\phi(\boldsymbol{x})$ and the symmetric localization of the points, it can be seen that by taking the same intensity 
$\Phi_D(t)$ for all non-centred Dirac delta points, reduces \eqref{eq:conditions for Phi} to two equations, for two unknown intensities, $\widetilde{\Phi}_D(t)$ and $\widetilde{\Phi}_C(t)$.

We will mainly discuss the case $n=1$ from this point onwards. Moreover, intuitively, if $\rho$ is very small, then the flux density function can be approximated well by a homogeneous flux density $\phi_0$. As our interest lies in the inhomogeneous flux density over the cell boundary, we select $\rho$ to be $1$ in this report.

The locations of the Dirac delta points are shown in Figure \ref{fig:1_2_Dirac}(a). The location of the non-centred Dirac delta point is $\boldsymbol{x}_1 = (x_c, y_c+r)$, where $(x_c, y_c)$ is the coordinate of the cell centre.  Solving \eqref{eq:expresion tilde Phi} -- \eqref{eq:conditions for Phi} under the conditions described above for $n=1$ leads to the following expressions for the intensities $\widetilde{\Phi}_D(t)$ and $\widetilde{\Phi}_C$:
\begin{equation}
    \label{Eq_dipole_sol}
    \left\{
    \begin{aligned}
        \widetilde{\Phi}_D(t) &= \frac{4\pi A(R+r)(R-r)}{(R+r)\exp\left\{-\frac{(R-r)^2}{4Dt}\right\}-(R-r)\exp\left\{-\frac{(R+r)^2}{4Dt}\right\}},\\
        \widetilde{\Phi}_C(t) &= 2\pi R\exp\left\{\frac{R^2}{4Dt}\right\}\left(\phi_0+A-\frac{2A(R+r)}{(R+r)-(R-r)\exp\left\{-\frac{Rr}{Dt}\right\}}\right). 
    \end{aligned}
    \right.
\end{equation}

\noindent Notice that, (\textit{i}) $r$ cannot be $0$, otherwise, there is no solution; (\textit{ii}) if $r = R$, then $\widetilde{\Phi}_D(t) = 0$, which indicates that the non-centred Dirac delta point is not needed at all; (\textit{iii}) at $t=0$, neither of the functions is integrable; (\textit{iv}) $\widetilde{\Phi}_D(t)$ is always non-negative, while $\widetilde{\Phi}_C(t)$ can be negative, in particular, if $\rho=1$ (i.e. $A = \phi_0$), then $\widetilde{\Phi}_C(t)$ is negative. In other words, the non-centred Dirac delta point is definitely a source point, while the centre point can be a sink or a source, which depends on the value of $\rho$. 

Substituting the solutions (Equation \eqref{Eq_dipole_sol}) into Equation \eqref{eq:expresion tilde Phi}, we managed to prove the convergence to the predefined sinusoidal flux density, which will \textit{not} be discussed here but in an upcoming article, in which more results will be shown.

\section{Numerical Results}
\label{sec:results}
In this section, we show numerical results of the quality of approximation of the spatial exclusion model with inhomogeneous prescribed flux density $\phi$ defined by Equation \eqref{Eq_phi} over the cell boundary of a single cell by means of (\textit{i}) a single-Dirac approach with source intensity $\Phi_0(t):=2\pi R\phi_0$ as in \cite{Peng2023} and (\textit{ii}) a multi-Dirac approach with symmetric localization at distance $r$ from the centre point, with intensities $\widetilde{\Phi}_C(t)$ at the centre and $\widetilde{\Phi}_D(t)$ at the non-centred points, given by \eqref{Eq_dipole_sol}. Recall, that determining the actual intensities $\Phi_i(t)$ was unfeasible at the moment. In the simulations we take fluctuations to be large: $\rho=A/\phi_0=1$. Furthermore, $\Omega=[-10,10]\times[-10,10]$, $\boldsymbol{x}_0=(0,0)$ and $R=1$, while $D=5$. For the initial condition, we assume there are no compounds in the extracellular environment, i.e. $u_0(\boldsymbol{x},0) = 0$ in the spatial exclusion model. Then, in the point source model, we define $\bar{u}_0(\boldsymbol{x},0) = 0$ as a spatial extension of $u_0(\boldsymbol{x},0)$ from $\Omega\setminus\Omega_C$ to $\Omega$. For non-zero initial conditions, more details are discussed in our previous work \cite{Peng2023}. Concerning the space and time discretization, FEM with Lagrange elements and backward Euler are utilized \cite{ vanKan2023,Vuik2023}. 

\begin{figure}[h!]
    \centering
    \subfigure[Flux density over the cell boundary with \newline $r = 0.05$ at various times $t$ $(n=1)$]{
    \includegraphics[width = 0.49\textwidth]{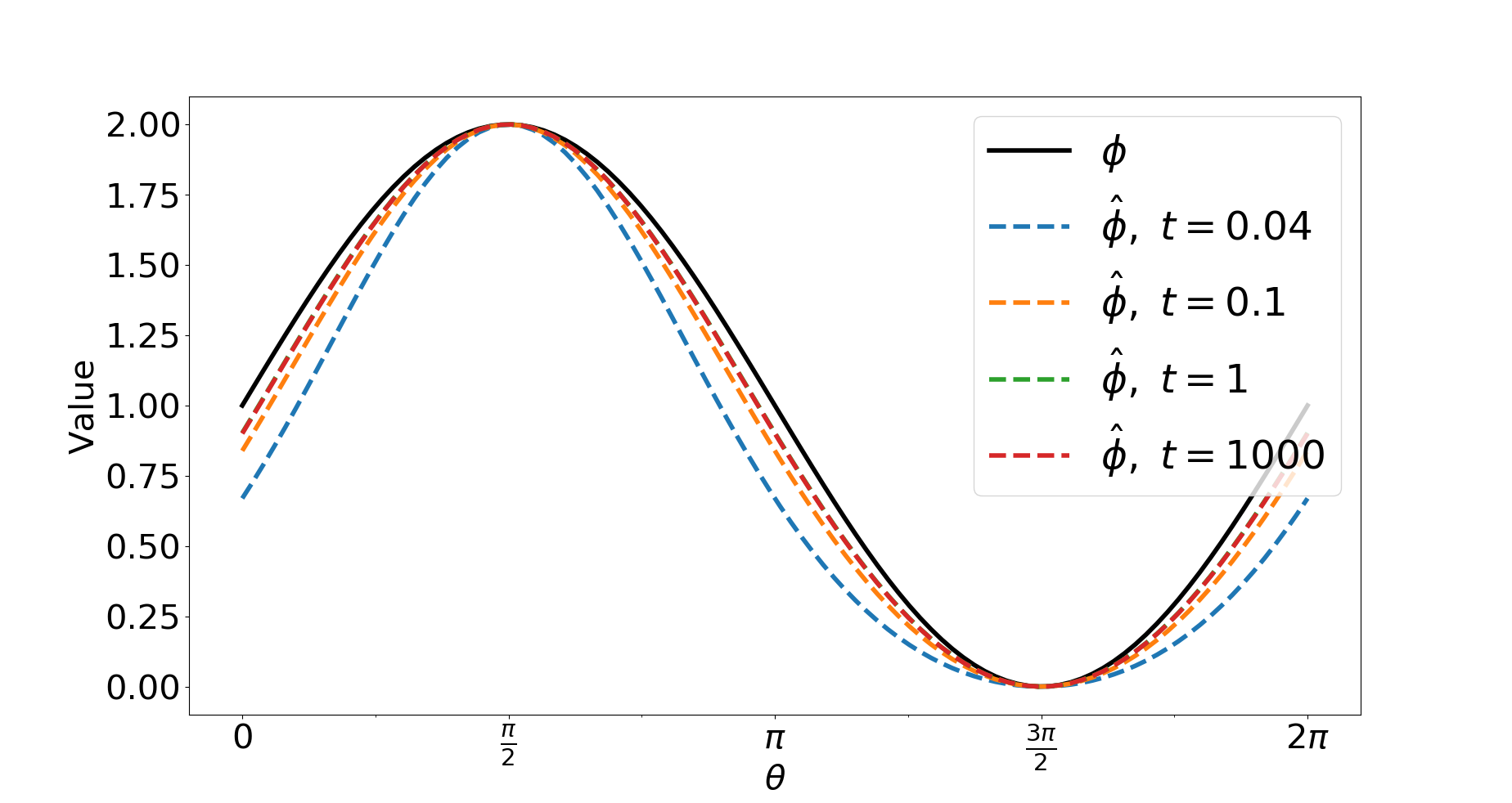}}
    \subfigure[Flux density over the cell boundary for various values of $r$ and $t\rightarrow+\infty$ $(n=1)$]{
    \includegraphics[width = 0.49\textwidth]{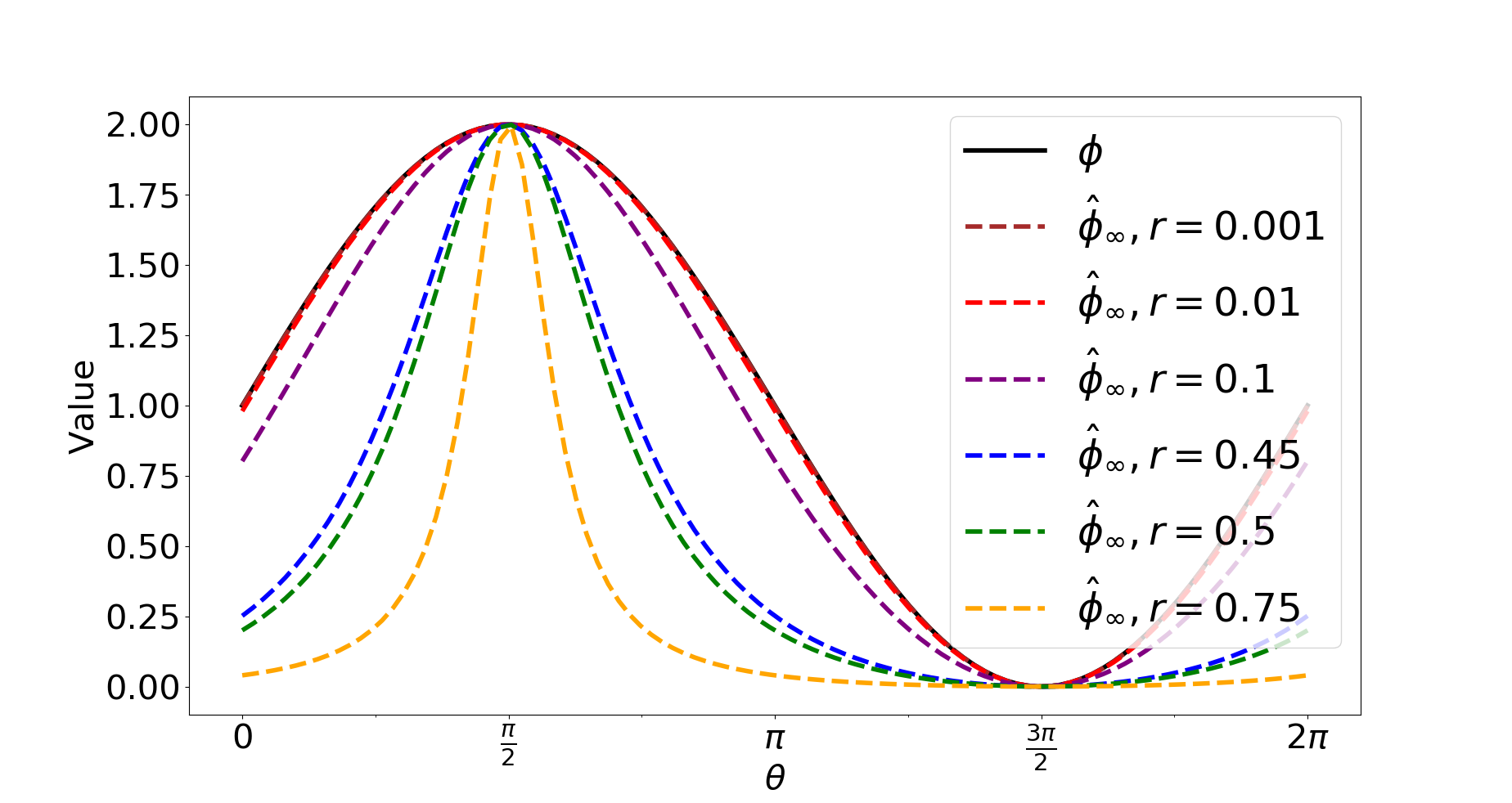}}
    \subfigure[Flux density over the cell boundary with \newline $r = 0.2$ at various times $t$ $(n=2)$]{
    \includegraphics[width = 0.49\textwidth]{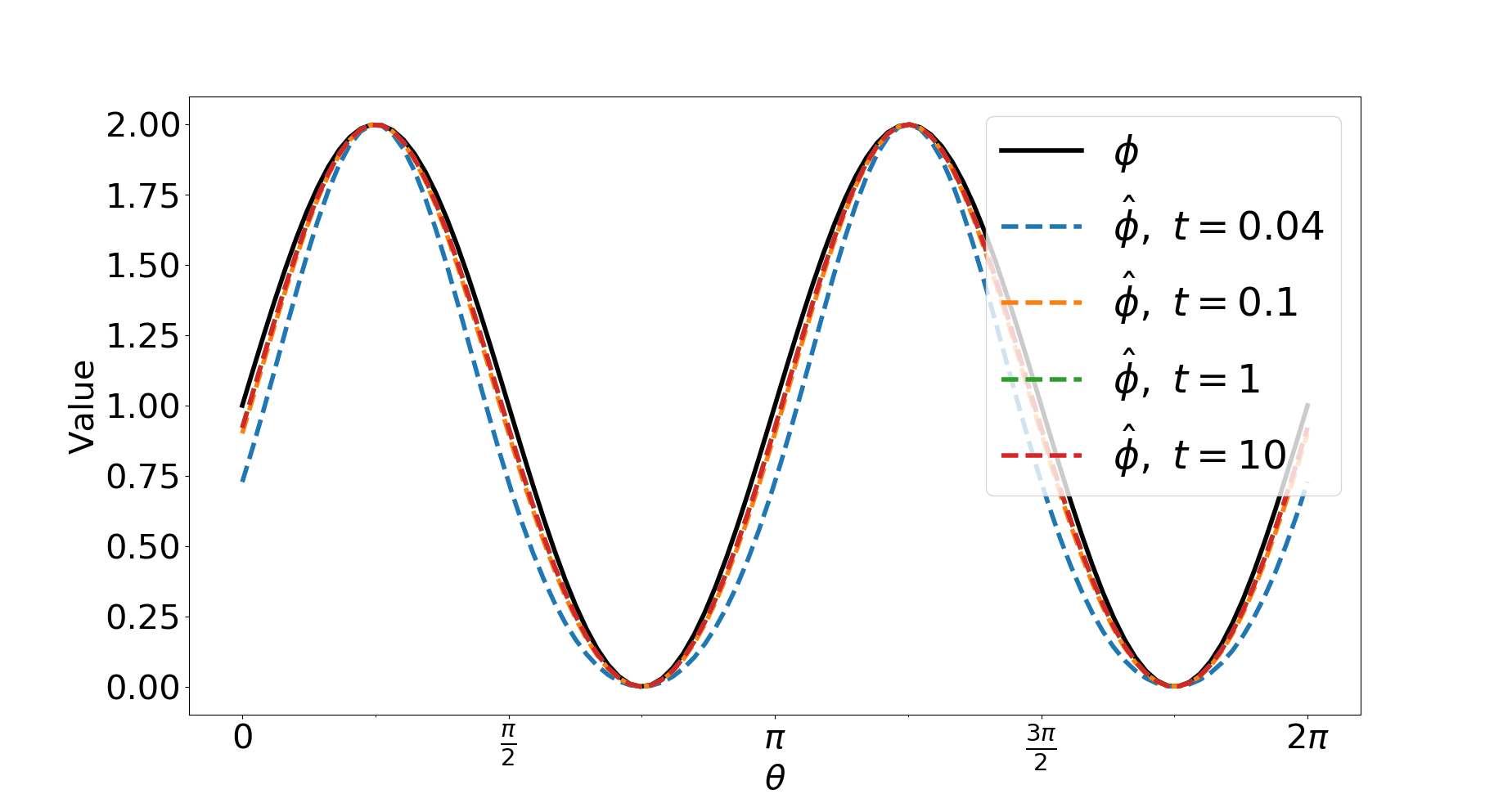}}
    \subfigure[Flux density over the cell boundary for various values of $r$ and $t\rightarrow+\infty$ $(n=2)$]{
    \includegraphics[width = 0.49\textwidth]{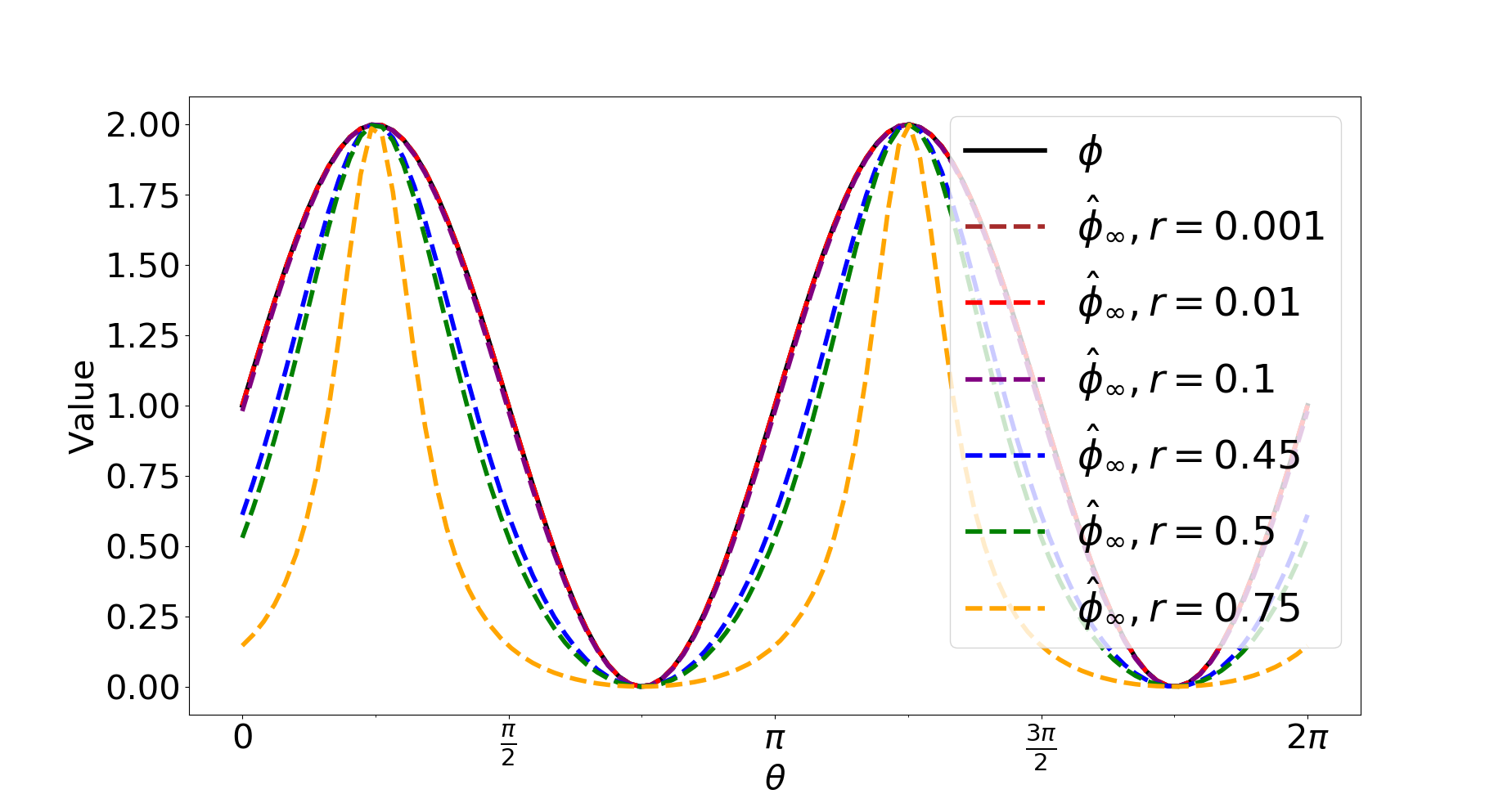}}   
    \caption{{\it Expression $\tilde{\phi}(\boldsymbol{x}_\theta,t)$ with symmetric configuration of the point sources and intensities $\widetilde{\Phi}_D(t)$ and $\widetilde{\Phi}_C(t)$ in comparison with $\phi(\boldsymbol{x}_\theta)$ given by \eqref{Eq_phi} with $n=1$. (a) The value of $r$ is fixed and $\tilde{\phi}(\boldsymbol{x}_\theta,t)$ is plotted at different $t$ for $n = 1$; (b) We take $t\rightarrow+\infty$ and the dashed curves show $\tilde{\phi}(\boldsymbol{x}_\theta,t)$ for large $t$ at different values of $r$ for $n = 1$. Panel (c) and (d) are similar plots as (a) and (b), respectively, with $n = 2$.}}
    \label{Fig_bnd_flux_n_1}
\end{figure}

In Figure \ref{Fig_bnd_flux_n_1} we show the flux density $\tilde{\phi}(\boldsymbol{x}_\theta,t)$ from this point source configuration in comparison to the objective $\phi(\boldsymbol{x}_\theta)$ for $n=1$ and $n=2$. Panel (a) and (c) shows the time evolution for fixed $r=0.01$, while Panel (b) and (d) shows the long-term behaviour for various $r$.

\noindent The numerical results suggest that when $r\rightarrow0^+$ and $t\rightarrow+\infty$,  $\tilde{\phi}(\boldsymbol{x}_\theta, t)$ converges to the prescribed $\phi(\boldsymbol{x})$. This convergence holds too for larger values of $n$; see Figure~\ref{Fig_bnd_flux_n_1}(c) and (d) where $n=2$ is shown. Convergence can be proven analytically as well, which will be shown in an upcoming publication. Moreover,
taking $r\rightarrow0^+$ and $t\rightarrow+\infty$ results in $\widetilde{\Phi}_D(t)\rightarrow+\infty$ and $\widetilde{\Phi}_C(t)\rightarrow-\infty$, then -- theoretically -- this multi-Dirac delta points setting is comparable to an electromagnetic dipole \cite{griffiths2005introduction}.

\begin{figure}[h!]
    \centering
    \subfigure[Error $\|u_S-u_P\|_{L^2(\Omega\setminus\Omega_C)}$]{
    \includegraphics[width = 0.48\textwidth]{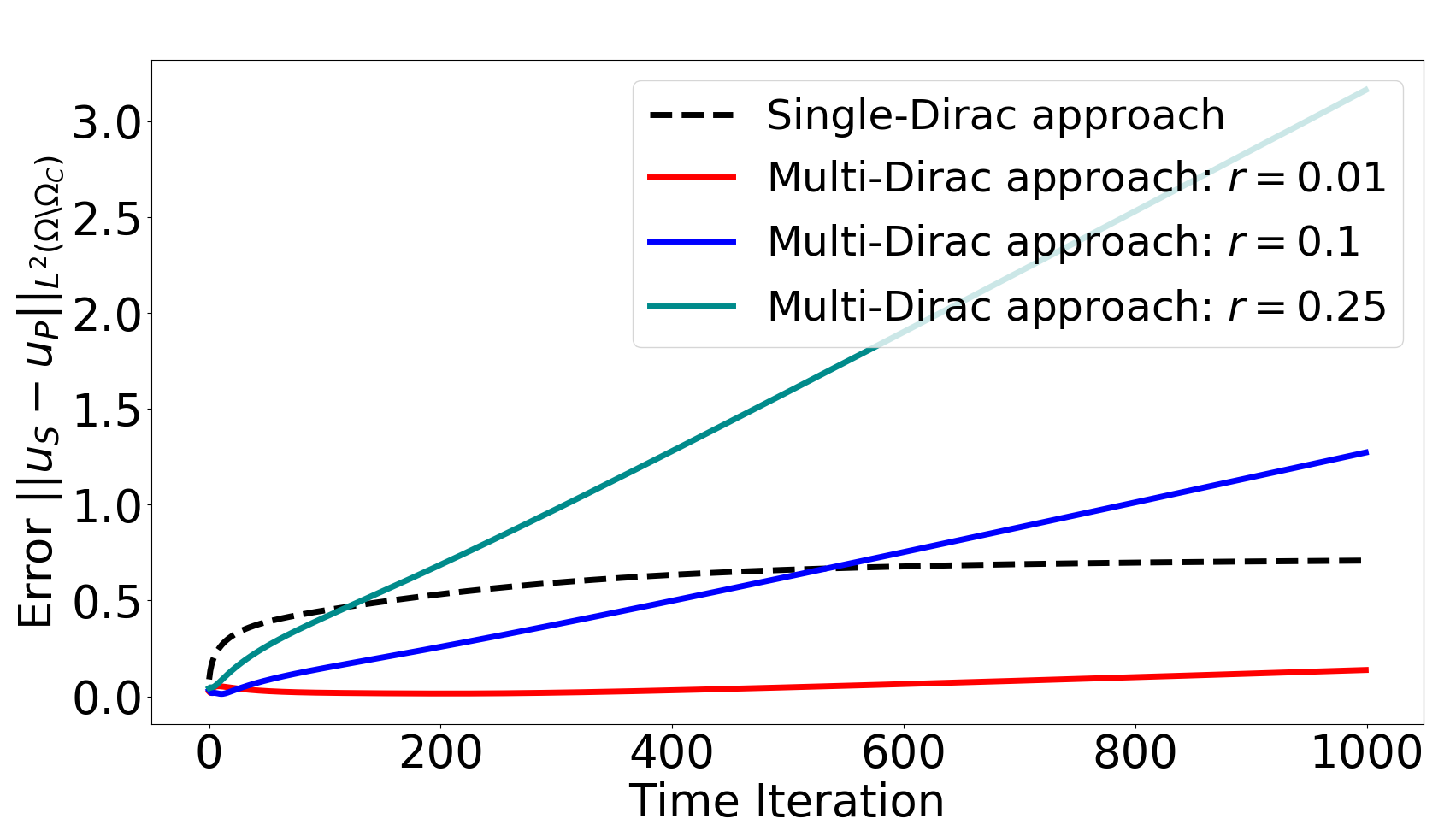}}
    \subfigure[Error $\|u_S-u_P\|_{H^1(\Omega\setminus\Omega_C)}$]{
    \includegraphics{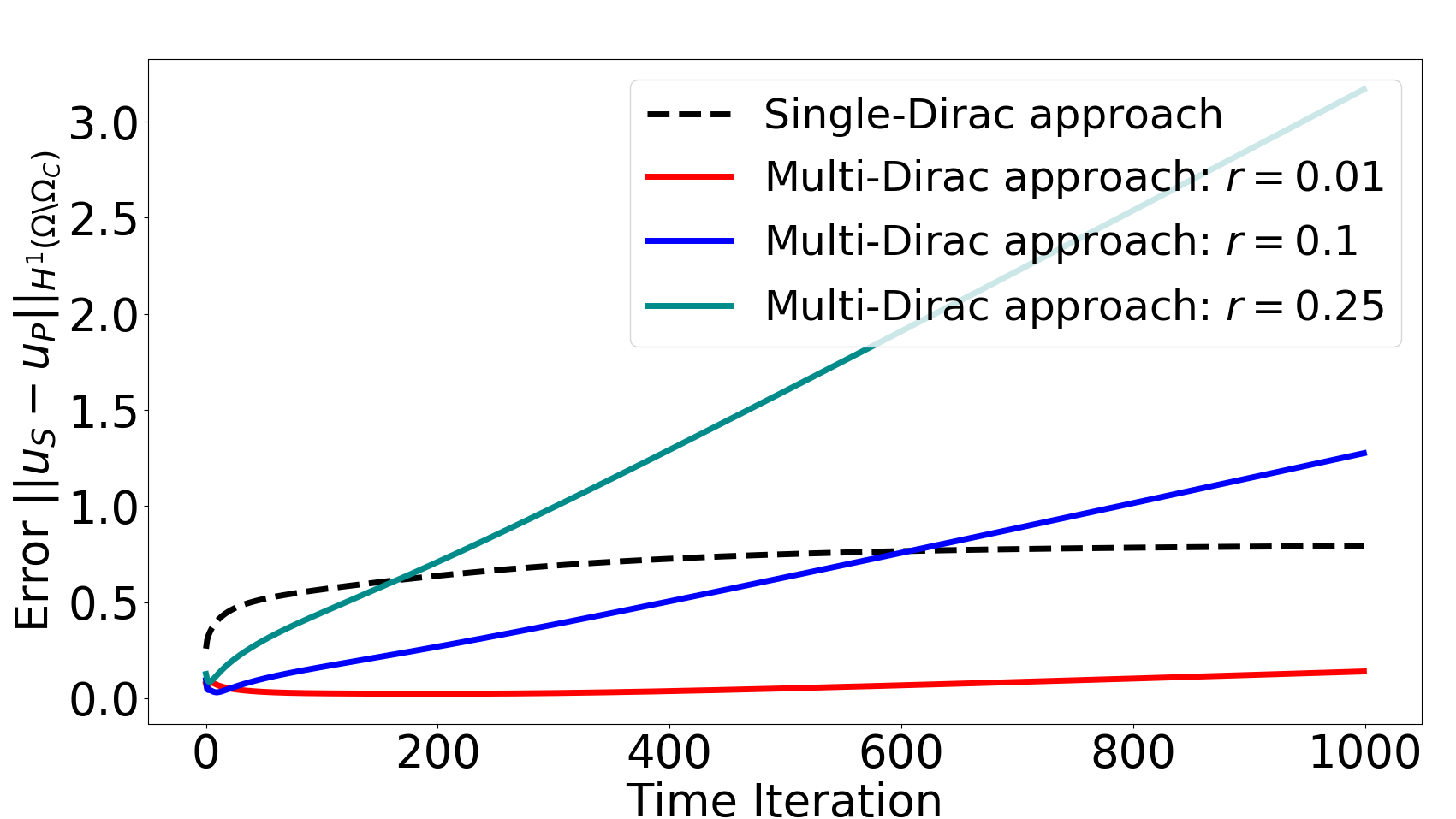}}
    \subfigure[$c^*(t)$]{
    \includegraphics[width = 0.48\textwidth]{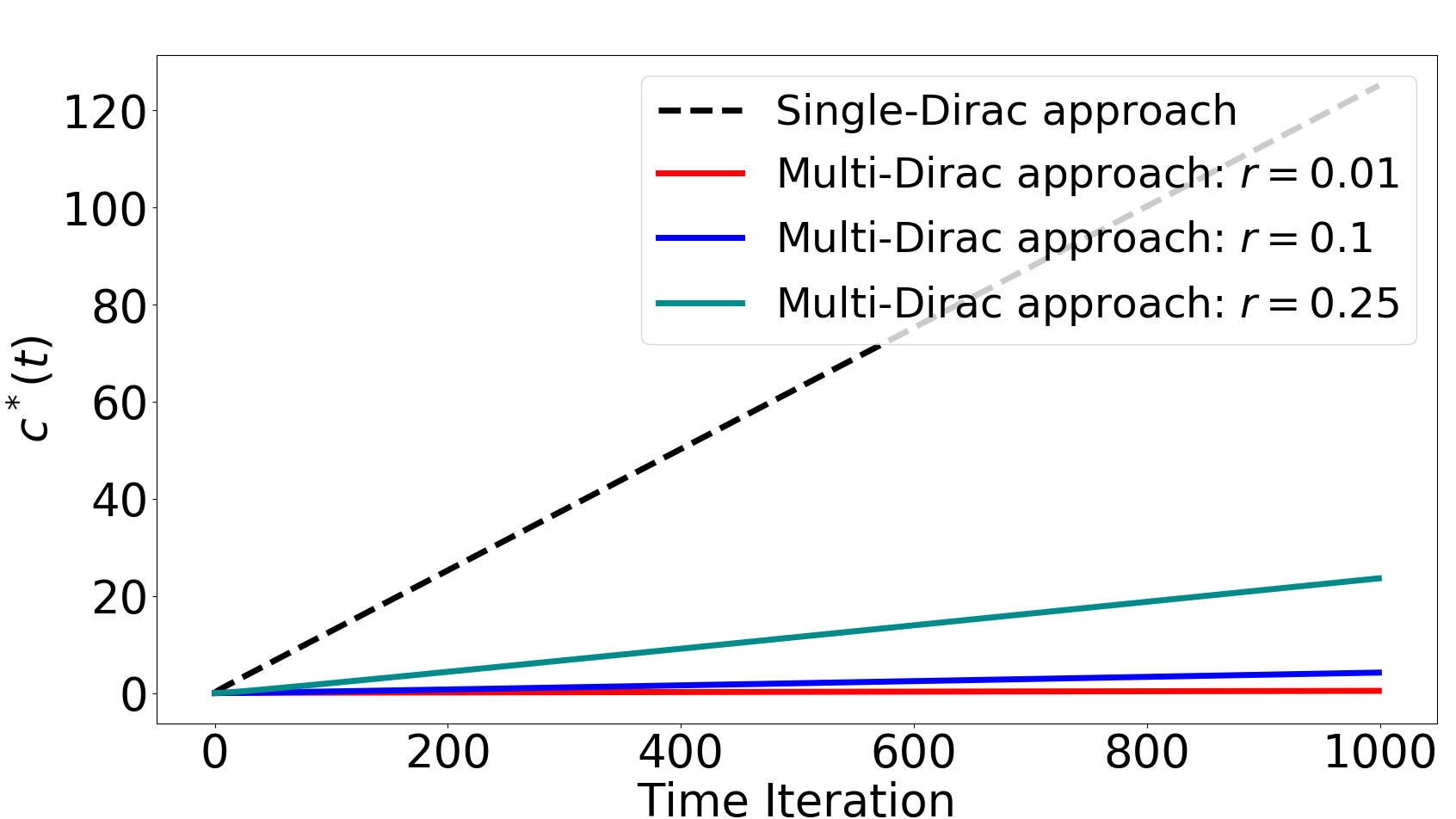}}
    \subfigure[$\|\phi(\cdot) - D\nabla u_P(\cdot,t    )\cdot\boldsymbol{n}\|_{L^2(\partial\Omega_C)}$]{
    \includegraphics[width = 0.48\textwidth]{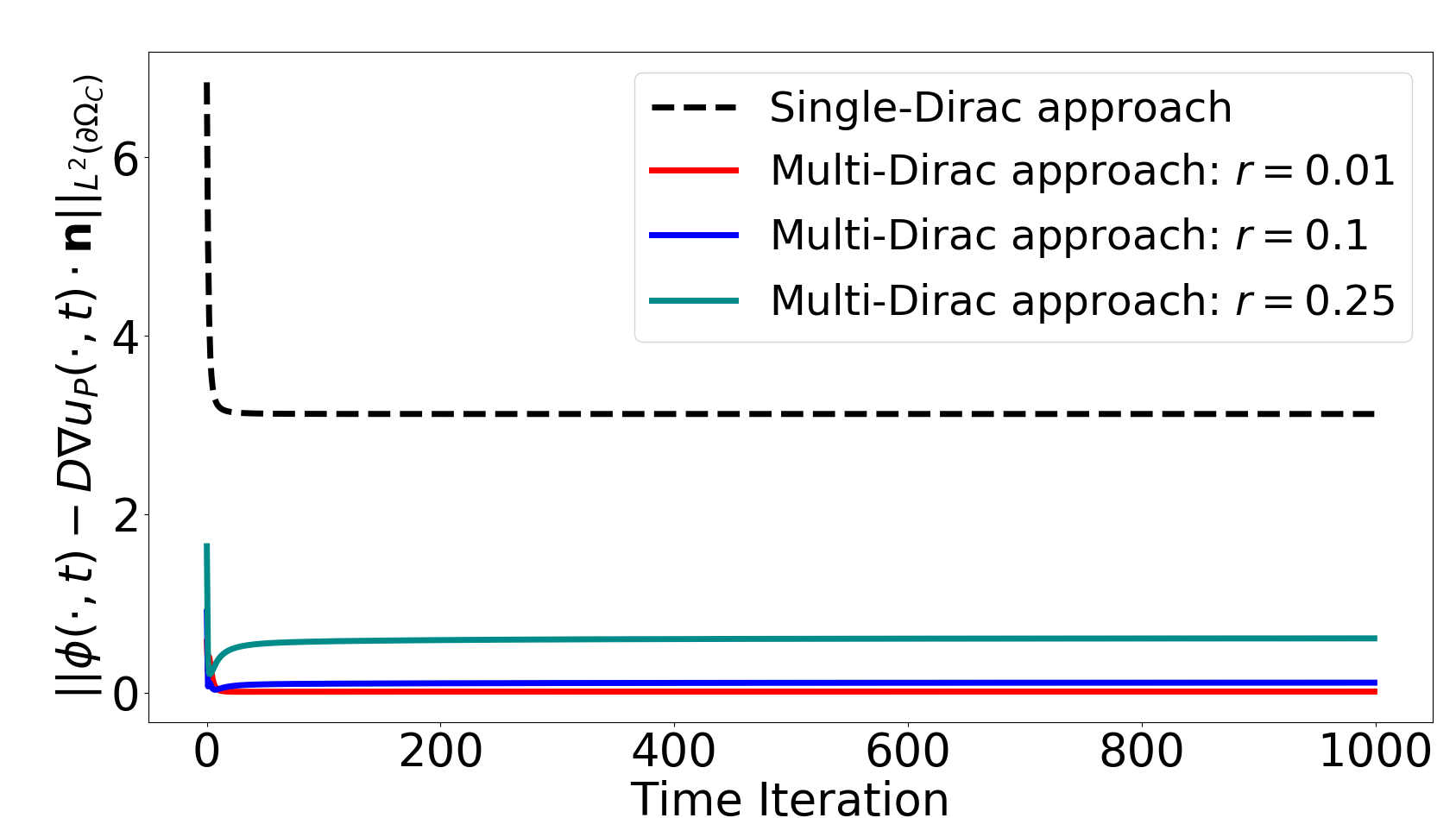}}
    \caption{{\it Multi-Dirac and single-Dirac approaches in the point source model with 'symmetric configuration' of points and $\varepsilon$-truncated intensities $\widetilde{\Phi}_D(t)$ and $\widetilde{\Phi}_c(t)$ with $n=1$ in comparison with the spatial exclusion model. The four quantities are plotted against the time iterations. They are computed by numerical integration after firstly interpolating the solutions to both models to the common mesh. $D=5$, timestep $\Delta t = 0.04$, terminal time $T = 40$, average mesh size $h = 0.0875$.}}
    \label{Fig_1_Dirac_ratio_1_D_5}
\end{figure}

 We shall now discuss evaluation of the effect of taking the approximate intensities $\widetilde{\Phi}_D(t)$ and $\widetilde{\Phi}_C(t)$ and which approach in the point source model can better approximate the spatial exclusion model. First, recall observation (\textit{iii}) in the preceding section, that both $\widetilde{\Phi}_D$ and $\widetilde{\Phi}_C$ are not integrable at $t=0$. Therefore, we truncate both functions for a small range of time $(0, \varepsilon)$ at the value at $\varepsilon$ in the numerical simulations to solve the boundary value problems in $(BVP_P)$. We take $\varepsilon = 0.01$. We consider the following quantitative measures for comparison: (\textit{i}) The $L^2-$norm difference $\|u_S-u_P\|_{L^2(\Omega\setminus\Omega_C)}$; (\textit{ii}) The $H^1-$norm difference $\|u_S-u_P\|_{H^1(\Omega\setminus\Omega_C)}$;\\ (\textit{iii}) The total  cumulative deviation in flux over the cell boundary, $c^*(t)$ from \cite{HMEvers2015}: 
\[
c^*(t) = \int_0^t\|\phi(\cdot) - D\nabla u_P(\cdot, s)\cdot\boldsymbol{n}\|_{L^2(\partial\Omega_C)}ds;
\]
(\textit{iv}) The flux deviation as a function of time  $\|\phi(\cdot) - D\nabla u_P(\cdot,t)\cdot\boldsymbol{n}\|_{L^2(\partial\Omega_C)}$.\\
The simulation results of these four quantities are shown in Figure \ref{Fig_1_Dirac_ratio_1_D_5}.

Generally speaking, the flux deviation over the entire cell boundary in the multi-Dirac approach is less than the one from single-Dirac approach,  see Figure \ref{Fig_1_Dirac_ratio_1_D_5}(d). It decreases with decreasing $r$, as expected from the improvement of flux approximation, shown in Figure \ref{Fig_bnd_flux_n_1}(b). However, if $r$ is not sufficiently small, the flux deviation at certain part of the circle might be very large, in particular at the beginning of the simulations, which can still result in a larger discrepancy between the solutions to the two models; see Figure \ref{Fig_1_Dirac_ratio_1_D_5}(a) and (b) when $r = 0.25$ and $r = 0.1$. In those cases, the single-Dirac approach outperforms the multi-Dirac for large time. Nevertheless, when $r$ is sufficiently small, e.g. $r = 0.01$, then the difference in the multi-Dirac approach is always less than using the single-Dirac approach in $t\in(0, 40)$.

\section{Conclusion and Discussion}
\label{sec:conclusion}
In this study, we investigated how to use the point source model to approximate the spatial exclusion model for a single cell when the predefined flux density over the cell boundary is inhomogenoeus, given by a specific sinusoidal function. This is a first explorative step towards settings in which there are more general flux inhomogeneities and more general cell shapes, and where there are many (and moving) cells. A point source model is expected to be computationally more efficient than a FEM approach, which requires frequent remeshing in the latter setting. Forthcoming publications on results in these directions are in preparation. 

We proposed using several Dirac-delta point sources to represent the cell, in a specific configuration and with well-chosen intensities, such that the inhomogeneous flux over the cell boundary can be approximated. To compute the intensity of the point sources, we enforced that the extreme points agree in the prescribed flux density function and in a modified flux expression for the point source model. We noticed numerically the convergence of this modified flux expression to the predefined flux density if $t\rightarrow+\infty$ and $r\rightarrow0^+$. Further numerical simulation showed that the multi-Dirac approach performs better in terms of solution difference metrics between the two models, if $r$ is sufficiently small. In summary, the multi-Dirac approach with a small $r$ is favoured when the compounds are secreted inhomogeneously via the cell boundary.

In forthcoming work, various research directions will be pursued. For instance, we already worked on multiple-cell situations, where we found that the distance between the cells will have a significant influence on the deviation between the two models \cite{Peng2023}; or one can think of the predefined flux density being spatial- and time-variant. To put in a nutshell, adding more complexity to the setup of the model automatically requires extra conditions to ensure appropriate consistency between these two diffusion models.

\bibliographystyle{abbrv}
\bibliography{authorsample}

% \section*{Appendix}
% \addcontentsline{toc}{section}{Appendix}

\end{document}